 
\baselineskip=17pt plus2pt minus1pt
\hsize=126truemm                     
\vsize=180truemm                      
\parindent=5truemm
\parskip=\smallskipamount
\mathsurround=1pt
\hoffset=2\baselineskip
\voffset=2truecm

%
%
\def\today{\ifcase\month\or
  January\or February\or March\or April\or May\or June\or
  July\or August\or September\or October\or November\or December\fi
  \space\number\day, \number\year}
%
 at 10truept

%
\newcount\dispno      
\dispno=1\relax       
\newcount\refno       
\refno=1\relax        
\newcount\citations   
\citations=0\relax    
\newcount\sectno      
\sectno=0\relax       
\newbox\boxscratch    
%

%
%
%
\def\Section#1#2{\global\advance\sectno by 1\relax%
\label{Section\noexpand~\the\sectno}{#2}%
\smallskip
\goodbreak
\setbox\boxscratch=\hbox{\bf Section \the\sectno.~}%
{\hangindent=\wd\boxscratch\hangafter=1
\noindent{\bf Section \the\sectno.~#1}\nobreak\smallskip\nobreak}}
%
\def\sqr#1#2{{\vcenter{\vbox{\hrule height.#2pt
              \hbox{\vrule width.#2pt height#1pt \kern#1pt
              \vrule width.#2pt}
              \hrule height.#2pt}}}}
\def\square{$\mathchoice\sqr34\sqr34\sqr{2.1}3\sqr{1.5}3$}
\def\endproof{~~\hfill\square\par\medbreak}
\def\noproof{~~\hfill\square}
%
%
\def\proc#1#2#3{{\hbox{${#3 \subseteq} \kern -#1cm _{#2 /}\hskip 0.05cm $}}}
\def\propcont{\mathchoice\proc{0.17}{\scriptscriptstyle}{}
                         \proc{0.17}{\scriptscriptstyle}{}
                         \proc{0.15}{\scriptscriptstyle}{\scriptstyle }
                         \proc{0.13}{\scriptscriptstyle}{\scriptscriptstyle}}
%

%
\def\normalin{\hbox{\raise0.045cm \hbox
                   {$\underline{\triangleleft }$}\hskip0.02cm}}
%
%
\def\'#1{\ifx#1i{\accent"13 \i}\else{\accent"13 #1}\fi}
%
%
%
\def\semidirect{\rlap{$\times$}\kern+7.2778pt \vrule height4.96333pt
width.5pt depth0pt\relax\;}
%
%
\def\prop#1#2{{\bf Proposition~\the\sectno.\the\dispno. }%
\label{Proposition\noexpand~\the\sectno.\the\dispno}{#1}\global\advance\dispno 
by 1{\it #2}\smallbreak}
\def\thm#1#2{{\bf Theorem~\the\sectno.\the\dispno. }%
\label{Theorem\noexpand~\the\sectno.\the\dispno}{#1}\global\advance\dispno
by 1{\it #2}\smallbreak}
\def\cor#1#2{{\bf Corollary~\the\sectno.\the\dispno. }%
\label{Corollary\noexpand~\the\sectno.\the\dispno}{#1}\global\advance\dispno by
1{\it #2}\smallbreak}
\def\defn{{\bf
Definition~\the\sectno.\the\dispno. }\global\advance\dispno by 1\relax}
\def\lemma#1#2{{\bf Lemma~\the\sectno.\the\dispno. }%
\label{Lemma\noexpand~\the\sectno.\the\dispno}{#1}\global\advance\dispno by
1{\it #2}\smallbreak}
\def\rmrk#1{{\bf Remark~\the\sectno.\the\dispno.}%
\label{Remark\noexpand~\the\sectno.\the\dispno}{#1}\global\advance\dispno
by 1\relax}
\def\proof{\noindent{\it Proof: }}
\def\numbeq#1{\the\sectno.\the\dispno\label{\the\sectno.\the\dispno}{#1}%
\global\advance\dispno by 1\relax}

\def\comm#1,#2{\left[#1{,}#2\right]}
\newdimen\boxitsep \boxitsep=0 true pt
\newdimen\boxith \boxith=.4 true pt 
\newdimen\boxitv \boxitv=.4 true pt
\gdef\boxit#1{\vbox{\hrule height\boxith
                    \hbox{\vrule width\boxitv\kern\boxitsep
                          \vbox{\kern\boxitsep#1\kern\boxitsep}%
                          \kern\boxitsep\vrule width\boxitv}
                    \hrule height\boxith}}
\def\square{\ \hbox{\vrule height7.5pt depth1.5pt width 6pt}\par}
\outer\def\square{\ifmmode\else\hfill\fi
   \setbox0=\hbox{} \wd0=6pt \ht0=7.5pt \dp0=1.5pt
   \raise-1.5pt\hbox{\boxit{\box0}\par}
}

\def\frac#1/#2{\leavevmode\kern.1em
              \raise.5ex\hbox{\the\scriptfont0 #1}\kern-.1em
              /\kern\.15em\lower.25ex\hbox{\the\scriptfont0 #2}}
\def\incnoteq{\lower.1ex \hbox{\rlap{\raise 1ex
     \hbox{$\scriptscriptstyle\subset$}}{$\scriptscriptstyle\not=$}}}
%
%
\def\mapright#1{\smash{
     \mathop{\longrightarrow}\limits^{#1}}}


\def\propcontup{\bigcup\!\!\!\rlap{\kern+.2pt$\backslash$}\,\kern+1pt\vert}
%
%
%
\def\label#1#2{\immediate\write\aux%
{\noexpand\def\expandafter\noexpand\csname#2\endcsname{#1}}}
%
\def\ifundefined#1{\expandafter\ifx\csname#1\endcsname\relax}
%
%
\def\ref#1{%
\ifundefined{#1}\message{! No ref. to #1;}%
 \else\csname #1\endcsname\fi}
%
%
\def\refer#1{%
\the\refno\label{\the\refno}{#1}%
\global\advance\refno by 1\relax}
%
%
\def\cite#1{%
\expandafter\gdef\csname x#1\endcsname{1}%
\global\advance\citations by 1\relax
\ifundefined{#1}\message{! No ref. to #1;}%
\else\csname #1\endcsname\fi}
%
%
 at 8truept      
%
%
%


\newread\aux
\immediate\openin\aux=\jobname.aux
\ifeof\aux \message{! No file \jobname.aux;}
\else \input \jobname.aux \immediate\closein\aux \fi
\newwrite\aux
\immediate\openout\aux=\jobname.aux

\headline={\ifnum\pageno<3{\hfill}\else{\ifodd\pageno\rhead\else\lhead\fi}\fi}
\def\lhead{A. Magidin\hfill}
\def\rhead{\hfill Dominions in decomposable varieties}
 
 at 8truept

\vfill 
\centerline{{\bf Dominions in decomposable varieties}\footnote*{The author was
supported in part by a fellowship from the Programa de Formaci\'on y
Superaci\'on del Personal Acad\'emico de la UNAM, administered by the
DGAPA.}}
\bigskip
{\obeylines 
Arturo Magidin 
Oficina 112 
Instituto de Matem\'aticas, UNAM 
Area de la Investigaci\'on Cient\'ifica 
Circuito Exterior 
Ciudad Universitaria 
04510 Mexico City, MEXICO 
e-mail: {\it magidin@matem.unam.mx}
}
\vfill\eject

\smallskip
{\parindent=20pt
\narrower\narrower
\noindent{Abstract. Dominions, in the sense of Isbell,
are investigated in the context of decomposable varieties of
groups. An upper and lower bound for dominions in such a variety is
given in terms of the two varietal factors, and the internal structure of the
group being analyzed. Finally, the following result is established: If
a variety ${\cal N}$ has instances of nontrivial
dominions, then for any proper subvariety ${\cal Q}$ of
${\cal G}roup$, ${\cal NQ}$ also has
instances of nontrivial~dominions.\par}}
\bigskip
\medskip


\bigskip\Section{Introduction}{intro}

Suppose that a group~$G$, a subgroup $H$ of~$G$, and a class ${\cal
C}$ of groups containing~$G$ are given. Are there any elements $g\in
G\setminus H$ such that any two morphisms between~$G$ and a ${\cal
C}$-group which agree on~$H$ must also agree on~$g$?

To put this question in a more general context, let~$\cal C$ be a full
subcategory of the category of all algebras (in the sense of Universal
Algebra) of a fixed type which is closed under passing to subalgebras.
Let $A\in {\cal C}$, and let~$B$ be a subalgebra of~$A$. Recall that,
in this situation, Isbell {\bf [\cite{isbellone}]} defines the
{\it
dominion of~$B$ in~$A$} (in the category ${\cal C}$) to be the
intersection of all equalizer subalgebras of~$A$
containing~$B$. Explicitly,
$${\rm dom}_A^{\cal C}(B)=\Bigl\{a\in A\bigm| \forall C\in {\cal C},\;
\forall f,g\colon A\to C,\ {\rm if}\ f|_B=g|_B{\rm\ then\ }
f(a)=g(a)\Bigr\}.$$
 
Therefore, the question with which we opened this discussion may be
rephrased in terms of the dominion of~$H$ in~$G$ in the category
of~context.

If $H={\rm dom}_G^{{\cal C}}(H)$ we say that the dominion of~$H$
in~$G$ (in the category~${\cal C}$) is {\it trivial}, and say it is
{\it nontrivial}~otherwise.
 
In this work we will study dominions when the category~${\cal C}$ is
a product of two proper nontrivial varieties of~groups. In
\ref{introandwreaths} we will recall the basic properties of varieties
that we will need; we refer the reader to Hanna Neumann's excellent
book {\bf [\cite{hneumann}]} for more information on varieties of
groups. Since wreath products are closely related to products of
varieties, and are used in the proofs of our results, we will also
recall some of their properties. In \ref{boundsfordoms} we will state
and prove the main results of this work, which give an upper and lower
bounds for the dominion of a subgroup in a product variety. Finally,
in \ref{applicationsnumber} we will use this result to prove that if a
variety~${\cal N}$ has instances of nontrivial dominions, so
will~${\cal NQ}$ for any variety ${\cal Q}\not={\cal G}roup$. Along
the way we will obtain some other results which are of interest in
their own~right.

The contents of this work are part of the author's doctoral
dissertation, which was conducted under the direction of Prof.~George
M.~Bergman, at the University of California at~Berkeley. It is my very
great pleasure to record and express my deep gratitude and
indebtedness to Prof.~Bergman. His advice and suggestions have
improved this work in ways too numerous to mention explicitly, and the
work itself would have been impossible without his help. He also
helped correct many mistakes; any errors that remain,
however, are my own~responsibility.

\Section{Preliminary results and wreath products}{introandwreaths}

All groups will be written multiplicatively, and all maps will be
assumed to be group morphisms, unless otherwise specified. Given a
group~$G$, the identity element of~$G$ will be denoted by~$e_G$, although
we will omit the subscript if it is understood from context. Given a
group $G$ and a subgroup~$H$, $N_G(H)$ denotes the normalizer of~$H$
in~$G$.

Recall that a variety of groups is a full subcategory of~${\cal
G}roup$ which is closed under taking arbitrary direct products,
quotients, and subgroups. Alternatively, it is the collection of all
groups (and all group morphisms between them) which satisfy a given set
of~identities.

We direct the reader to~{\bf [\cite{nildomsprelim}]} for the basic
properties of dominions. I mention the most important of them, which
are not hard to verify: ${\rm dom}_G^{{\cal C}}(-)$ is a closure
operator on the lattice of subgroups of~$G$; if ${\cal C}$ is closed
under quotients, then normal subgroups are dominion-closed;
the dominion construction respects finite direct products; and
the dominions construction respects quotients in a variety. That is,
if ${\cal V}$ is a variety of groups, $G\in {\cal V}$, and $H$ is a
subgroup of~$G$, $N$ a normal subgroup of~$G$ contained in~$H$, then
$${\rm dom}_{G/N}^{{\cal V}}(H/N) = {\rm dom}_G^{\cal V}(H)\bigm/ N.$$

Given a variety~${\cal V}$, a group~$G\in {\cal V}$ and a subgroup $H$
of~$G$, the {\it amalgamated coproduct (in ${\cal V}$) of $G$ with
itself over~$H$}, denoted by~$G\amalg_H^{\cal V} G$ is the universal
${\cal V}$-group, equipped with embeddings $\lambda,\rho\colon G\to
G\amalg_H^{\cal V} G$ such that $\lambda|_H =\rho|_H$.  If $H$ is the
trivial subgroup, the resulting object is the usual coproduct
in~${\cal V}$, which we denote simply by $G\amalg^{\cal V}G$. It is
not hard to verify that in a variety, the dominion of $H$ is 
the equalizer of the two canonical embeddings into the amalgamated
coproduct.

Given two varieties of groups, $\cal N$ and $\cal Q$, recall that the
{\it product variety ${\cal V} = {\cal NQ}$} is the variety of all
groups which are extensions of a group in~${\cal N}$ by a group
in~${\cal Q}$; that is, it consists of all groups~$G$ which have a
normal subgroup $N\in {\cal N}$, such that $G/N\in{\cal Q}$. 

We will say that a variety ${\cal V}$ is {\it nontrivial} iff ${\cal
V}\not={\cal G}$ and ${\cal V}\not={\cal E}$, and we will call it {\it
trivial} otherwise. A variety ${\cal V}$ {\it factors nontrivially}
(or {\it is decomposable}) if it can be expressed as the product of
two nontrivial varieties.

The semigroup of varieties of groups has the structure of a
cancellation semigroup with 0 and 1. The zero element is the variety
of all groups (denoted by ${\cal G}$), while the identity is the
trivial variety, consisting only of the trivial group and denoted
by~${\cal E}$. Furthermore, every variety other than ${\cal G}$ can be
uniquely factored as a product of a finite number of indecomposable
varieties (in a unique order), with ${\cal E}$ having the empty
factorization, so that the semigroup with neutral element of varieties
other than~${\cal G}$ is freely generated by the indecomposable
varieties.  See Theorems~21.72, 23.32 and~23.4 in {\bf
[\cite{hneumann}]}.

Given a variety ${\cal V}$ and a group $G$ (not necessarily in~${\cal
V}$), we will denote by~${\cal V}(G)$ the verbal subgroup of~$G$
associated to~${\cal V}$. This is the subgroup generated by all
values of the words ${\bf w}$ which are laws of~${\cal V}$. In
particular, $G\in {\cal V}$ if and only if~${\cal V}(G) = \{e\}$. We
also note the universal property associated to ${\cal V}(G)$: for any
normal subgroup $N\triangleleft G$, $G/N\in{\cal V}$ if and only
if~${\cal V}(G)\subseteq N$.

The term {\it $n$-generator group} will mean that the
group in question can be generated by $n$ elements, but may in fact
need less. The term {\it $n$-variable word} will refer to a word in
$x_1,\ldots,x_n$. This convention relies on the fact that in a law
involving $n$ variables, the name of these variables is~immaterial.

We also recall the definition of the {\it wreath product} of two
groups. Given groups $N$ and $K$, the {\it regular (unrestricted)
wreath product of~$N$ by~$K$}, denoted by $N\wr K$, is constructed as
follows:

We take $N^K$, the cartesian power of~$N$ consisting of all functions
$\phi\colon K\to N$ (not necessarily group morphisms) multiplied
componentwise. For each element $k\in K$, we let $\beta_k\colon N^K\to
N^K$ be the mapping that takes $\phi\in N^K$ to~$\phi^k$, by
$$\phi^k(y) = \phi(yk^{-1}) \qquad {\rm for\ all\ }y\in K.$$

Then $\beta_k$ is an automorphism of $N^K$, and (if $N\not=\{e\}$) the
set of these automorphisms is a group isomorphic to~$K$. Let~$P$ be
the semidirect product of $N^K$ by this group of automorphisms: that
is, we take pairs $(k,\phi)$, with $k\in K$ and $\phi\in N^K$, with
multiplication
$$(k,\phi) (\ell,\psi) = (k\ell,\phi^{\ell}\psi),$$ and identify the
group of elements $(k,e)$ with~$K$, and the group
of elements $(e,\phi)$ with $N^K$ so that elements
of~$P$ become products $k\phi$, and we have $\phi^k=k^{-1}\phi k$. In
this notation, $P$ is also called the {\it complete} (or {\it
unrestricted}) wreath product.

If we replace $N^K$ with $N^{(K)}$, the subgroup of
functions $\phi\colon K\to N$ such that $\phi(k)=e$ for almost all
$k\in K$, the resulting group is called the {\it
restricted wreath product}, and we shall denote it by $N\>{\rm
wr}\>K$. Clearly, $N\>{\rm wr}\>K$ is a subgroup of~$N\wr K$.

In general, if $\Omega$ is a $K$-set, we may form the product of
$|\Omega|$ copies of~$N$, and then let~$K$ act on this product by
permuting the components. This is called the $\Omega$-wreath product
of~$N$ by~$K$, and is denoted by~$N\wr_{\Omega} K$. We see that the
regular wreath product is nothing more than a special case of this
situation, where $\Omega$ is $K$ itself, with the action being given
by the right regular~action.

Next we recall some basic results about product varieties and wreath
products. Their proofs are straightforward and therefore we omit them,
although the reader will be able to provide them herself if she
so~desires.

\prop{verbproduct}{{\rm (Proposition~21.12 in~{\bf [\cite{hneumann}]})} Let
${\cal N}$ and ${\cal Q}$ be two varieties of 
groups. Given a group~$G$, the verbal subgroup corresponding to the
variety ${\cal NQ}$ is given by ${\cal N}({\cal Q}(G))$.\noproof}

\lemma{ngengroup}{{\rm (Lemma~16.1 in~{\bf [\cite{hneumann}]})} An
$n$-generator group~$G$ belongs to the variety~${\cal V}$ if and only
if it satisfies the $n$-variable laws of~${\cal V}$.\noproof}

We let ${\cal V}^{(n)}$ denote the variety defined by the set of
$n$-variable laws of~${\cal V}$.

\thm{nvarlaws}{{\rm (Theorems~16.21 and 16.23 in~{\bf
[\cite{hneumann}]})} The variety ${\cal V}^{(n)}$ consists of all
groups whose $n$-generator subgroups belong to~${\cal V}$. Therefore,
${\cal V}^{(n)}\supseteq {\cal V}$, and ${\cal V}^{(n)}={\cal V}$ if
and only if ${\cal V}$ can be defined by $n$-variable laws. Also
$${\cal V}^{(1)}\supseteq {\cal V}^{(2)}\supseteq \cdots \supseteq
{\cal V}^{(n)}\supseteq \cdots \supseteq {\cal V} = \bigcap_{i\geq 1}
{\cal V}^{(i)}.\eqno\noproof$$}

\thm{wreaths}{Let $N\in{\cal N}$ and $Q\in{\cal Q}$. Then $N\wr Q\in
{\cal NQ}$ and $N\>{\rm wr}\>Q\in{\cal NQ}$.}

\proof $N\wr Q$ is an extension of $N^Q$ by~$Q$. $N^Q$ is an ${\cal
N}$-group, since it is the direct product of copies of $N$, and $N\in
{\cal N}$. Therefore, $N\wr Q$ is the extension of an ${\cal N}$-group
by a~${\cal Q}$-group, hence an element of ${\cal NQ}$.

Finally, since $N\>{\rm wr}\>Q$ is a subgroup of $N\wr Q$, it is also
an element of ${\cal NQ}$.\endproof

Let $N$ and~$K$ be two groups, and let $f\colon N\to M$ be a group
morphism. Given an element $k\phi\in N\wr K$, we may obtain an element
of $M\wr K$ induced by~$f$, namely $k(f\circ\phi)$; note that $\phi$
is a map from $K$ to~$N$, so $f\circ\phi$ is a map from $K$ to~$M$,
that is an element of~$M^K$. We let the induced map from $N\wr K$ to
$M\wr K$ be denoted by~$f^*$. With a straightforward verification we~obtain:

\thm{propsofwreaths}{{\rm (See 22.21--22.14 in~{\bf [\cite{hneumann}]})} Let
$N$ and~$K$ be two groups, and let $f\colon N\to M$ be a group
morphism. Then the induced map $f^*\colon N\wr K
\to M\wr K$ is a group homomorphism. Furthermore, if $f$ is an
embedding, then so is~$f^*$.\noproof}

\thm{extensions}{{\rm (Kaloujnine and Krasner {\bf [\cite{wreathext}]})} The
complete wreath product $A\wr B$ contains an isomorphic copy of every
group that is an extension of~$A$ by~$B$.}

\proof Let
$$1\,\mapright{}\,A\,\mapright{\alpha}\,G\,\mapright{\pi}\,B\,\mapright{}\,1$$
be an exact sequence of groups, and let $\tau$ be a transversal
of~$\pi$ (that is, $\pi\circ\tau = {\rm id}_B$; note that $\tau$ is
not necessarily a group morphism). Let $T$ be the image
of~$\tau$. Thus, $T$ is a set of coset representatives for $\alpha(A)$ in~$G$.
We now define a mapping $\gamma\colon
G\to A\wr B$ as follows: for $g\in G$, let
$$\eqalign{\gamma(g) &= \pi(g)\varphi_g,\qquad\hbox{where for all
$y\in B$}\cr
\varphi_g(y)&= \alpha^{-1} \Bigl( \tau\bigl(y\pi(g)^{-1}\bigr) g
\tau(y)^{-1}\Bigr);\cr}$$
one can check that the right hand side is meaningful. A
straightforward computation will confirm that $\gamma$ is a
monomorphism.\endproof

\rmrk{dependsonT} Note that the embedding $\gamma$ obtained in
\ref{extensions} depends on the choice of transversal~$\tau$ (or
equivalently, on the choice of coset representantives~$T$).

\rmrk{categoricaladdenda} Given a fixed group~$B$, we can make the set
of groups with maps to~$B$ into a category, using for morphisms the
group homomorphisms that make commuting triangles with the maps
into~$B$. Given $G$ as an extension of~$A$ by~$B$, we have a given map
$\pi\colon G\to B$, and we also have a canonical map from $A\wr B$
to~$B$. It is easy to verify that the embedding of~$G$ into $A\wr B$
given in \ref{extensions} is actually a morphism in this category.

\Section{Upper and lower bounds for the dominion in a product
variety}{boundsfordoms}

We now turn our attention to the behavior of dominions in product
varieties. Throughout this section, we will work inside a variety
of groups ${\cal V}$ which can be decomposed as a product ${\cal
V}={\cal NQ}$ of two nontrivial varieties ${\cal N}$ and~${\cal Q}$,
unless otherwise specified.

\thm{secondfactor}{{\rm (See {\bf [\cite{mckay}]})}
Suppose that ${\cal V}={\cal NQ}$, where
$\cal N$ is a nontrivial variety. Let $G\in \cal Q$, and $H$ a
subgroup of~$G$. Then ${\rm dom}^{\cal V}_{G}(H)= H$.}

\proof Let $G$ and $H$ be as in the statement. 
The key idea is that if one can find a group~$N\in{\cal N}$ and an
action of~$G$ on~$N$ such that the stabilizer in~$G$ of some $n\in N$
is exactly~$H$, then we may form the semidirect product $N\semidirect
G$ determined by this action, which will lie in~${\cal V}$. In this
semidirect product, $H$ becomes the centralizer in~$G$ of the
element~$n$ of the semidirect product. Then the equalizer of the
inclusion of~$G$ into the semidirect product, and the composite of
this inclusion with conjugation by~$n$, is precisely~$H$, which would
prove that~$H$ is equal to its own dominion in~$G$ in the
variety~${\cal V}$.

We can obtain such an ${\cal N}$-group by letting~$N$ be a direct
product of~$|G/H|$ copies of a single nontrivial group~$M\in{\cal N}$,
with the action 
of~$G$ permuting the copies by left multiplication on the left cosets
of~$H$. We then let~$n$ be an
element of $N\subseteq M\wr_{G/H} G$ with support equal
to~$\{H\}$, which finishes the~proof.\endproof

What about groups in ${\cal V}$ but not ${\cal Q}$? We have a
partial answer in the following three results:

\thm{mckaygeneralized}{Let ${\cal V}={\cal NQ}$ be a variety, where
$\cal N$ is a nontrivial variety. Let $G\in{\cal
V}$, and let
$N\triangleleft G$, with $N\in{\cal N}$ and ${G/N\in {\cal Q}}$. Then
for
all subgroups~$H$ of~$G$, ${\rm dom}^{\cal V}_{G}(H)\subseteq NH$. In
particular, 
$${\rm dom}_G^{\cal V}(H)\subseteq {\cal Q}(G)H.$$}

\proof Since $H\subseteq NH$, we have 
that ${\rm dom}_{G}^{\cal V}({H})\subseteq{\rm dom}_{G}^{\cal
V}({NH})$. Also, $$NH/N\subseteq G/N\in {\cal Q}.$$

By \ref{secondfactor}, ${\rm dom}_{G/N}^{\cal V}({NH/N}) =
NH/N$. Since dominions respect quotients, we
know that ${\rm dom}_{G/N}^{\cal V}({NH/N}) = ({\rm
dom}_{G}^{\cal V}({NH}))/N$. Therefore,
$${\rm dom}_{G}^{\cal V}({H})\subseteq {\rm dom}_{G}^{\cal V}({NH}) =
NH$$
completing the~proof.\endproof

\lemma{bighalf}{Let ${\cal V}={\cal NQ}$ be a nontrivial factorization
of~${\cal V}$, and let $G\in {\cal V}$, $H$ a subgroup of~$G$, and
$N={\cal Q}(G)$. Let $D={\rm dom}_N^{\cal N}(N\cap H)$. Then
$$HD=\langle H,D\rangle \subseteq {\rm dom}_G^{\cal NQ}(H).$$}

\proof Let $K\in {\cal V}$ and let $f,g\colon G\to K$ be two group
morphisms, such that $f|_H = g|_H$.

Since $N={\cal Q}(G)$, we must have $f(N),g(N)\subseteq {\cal Q}(K)$. In
particular, $f|_N$ and $g|_N$ are actually maps of $\cal N$-groups. We
also have $f|_{H\cap N} = g|_{H\cap N}$, so we must have that $f|_N$
and~$g|_N$ agree on ${\rm dom}_N^{\cal N}(N\cap H)=D$. In particular,
$f$ and $g$ agree on~$D$. Therefore,
$$\langle H,D\rangle \subseteq {\rm dom}_G^{\cal NQ}(H).$$

Next we claim that $\langle H,D\rangle = HD$. Indeed, $H\cap N$ is
normal in~$H$, and since
the dominion certainly respects isomorphisms of group
and subgroup pairs, the ${\cal N}$-dominion $D$ of $H\cap N$ in~$N$ is
likewise fixed by the automorphisms of~$N$ induced by elements
of~$H$. Thus $H\subseteq N_G(D)$, so $\langle H,D\rangle = HD$,
as~claimed.\endproof

With slightly more information, we can obtain a stronger result:

\thm{bigone}{Let ${\cal V}={\cal NQ}$ be a nontrivial factorization
of~${\cal V}$, and let $G\in {\cal V}$, $H$ a subgroup of~$G$, and
$N={\cal Q}(G)$. Let $D={\rm dom}_N^{\cal N}(N\cap H)$.
If $N_G(D)N=G$, then
$${\rm dom}_G^{\cal NQ}(H) =  HD.$$}

\proof First note that $\langle H,D\rangle=HD$ as in the proof of
\ref{bighalf}.

By \ref{bighalf}, we already know that $$HD\subseteq{\rm dom}_G^{\cal
NQ}(H).$$
To prove that under the added hypothesis that
$N_G(D)N=G$ we also get the inclusion ${\rm dom}_G^{\cal
NQ}(H)\subseteq HD$, we need to define a transversal of~$N$ in~$G$. We
claim there is a transversal $\tau\colon G/N\to G$ with the following
properties:

{\parindent=45pt
\item{(\numbeq{proptransone})} $\tau (N) = e$.\par
\item{(\numbeq{proptranstwo})} $\tau(yN) \in N_G(D)$ for all $yN\in
G/N$.\par
\item{(\numbeq{proptransthree})}{For every $h\in H,\  y\in G$,
there exists an element $h'\in H$ such that $\tau(yh^{-1}N) =
\tau(yN)h'^{-1}$.}\par}

Since $N_G(D)N=G$, conditions (\ref{proptransone}) and
(\ref{proptranstwo}) are easy to achieve. The reason for
(\ref{proptransthree}) will become apparent during the course of the proof.

We now establish the existence of a transversal satisfying
(\ref{proptransone})--(\ref{proptransthree}). 
Let us look at the left action of~$H$ on the set of
cosets of~$N$, under which $h\in H$ takes the coset $tN$ to the coset
$tNh^{-1} = th^{-1}N$. Since $N$ is normal, this is a well defined
action. Thus, the set of cosets of~$N$ is divided into orbits by
the~$H$-action.

For each $H$-orbit, we
first define $\tau$ to take some arbitrary coset $tN$ in that orbit to
a representative in $N_G(D)$, which we now choose once and for all,
making sure to select $e$ as a representative for~$N$ to satisfy
(\ref{proptransone}).  For any other coset $t'N$ in the same orbit,
there exists an element $h\in H$ such that $t'\equiv
\tau(tN)h^{-1}
\pmod{N}$, because this corresponds to the $H$-action on the set of
cosets. Choose such an $h$ for each coset (the choice of~$h$ is only
determined up to congruence modulo~\hbox{$H\cap N$}), and define
$\tau(t'N) = \tau(tN)h^{-1}$.  Since $H\subseteq N_G(D)$, this
guarantees that both (\ref{proptranstwo}) and (\ref{proptransthree})
are satisfied by~$\tau$. Note that the above choices require the use
of the Axiom of Choice if the number of orbits is not~finite.

Let $\pi\colon G\to G/N$ be the
canonical projection onto the quotient. For simplicity, we will now
write the cosets using their chosen representatives; that is, whenever
we write such a coset as $tN$ it will be understood that~$t$ is the
chosen representative of that coset. If we wish to represent the coset
of an arbitrary element $y\in G$, where $y$ is not the chosen
representative of its coset, we will call this~$\pi(y)$.

Since $G$ is an extension of~$N$ by $G/N$, we have an embedding defined
as in \ref{extensions} (and depending on our choice of~$\tau$),
$\gamma\colon G\to N\wr G/N$; since $N\in {\cal N}$, and $N={\cal
Q}(G)$, we have $G/N\in {\cal Q}$, so~$N\wr G/N \in {\cal NQ}$.

Consider the group $N\amalg_D^{\cal N} N$.
Let $\lambda,\rho\colon N\to N\amalg_D^{\cal N} N$ be
the two immersions of $N$ into $N\amalg_D^{\cal N} N$. Since $D$ is
its own dominion in~$N$ (in the variety~${\cal N}$), 
we have $\lambda(n)=\rho(n)$ if and only if $n\in D$.

By \ref{propsofwreaths}, the maps $\lambda$ and~$\rho$ induce two maps
$$\lambda^*,\rho^*\colon N\wr (G/N) \to (N\amalg_D^{\cal N} N)\wr
(G/N).$$ 
We now consider the composite maps $\lambda^*\circ\gamma$ and
$\rho^*\circ\gamma$. 

Let $n\in N$. By definition of~$\gamma$, we have
$\gamma(n)=\varphi_n$, where $\varphi_n\colon G/N\to N$ is given~by
$$\eqalign{\varphi_n(yN) &= \tau\left( yN\pi(n)^{-1}\right) n \tau(yN)^{-1}\cr
&=\tau(yN)n\tau(yN)^{-1}\cr
&=yny^{-1}.\cr}$$

Therefore, $\lambda^*\circ\varphi_n (yN) = \lambda(yny^{-1})$, and
$\rho^*\circ\varphi_n (yN) = \rho (yny^{-1})$. Comparing the values at
$y=e$, we see that $\lambda^*\circ\gamma(n)=\rho^*\circ\gamma(n)$ only
if $n\in D$. We claim that in fact $\lambda^*\circ\gamma(n) =
\rho^*\circ\gamma(n)$ if and only if $n\in D$.

To establish this claim, let $d\in D$, and let $yN\in G/N$. Then we
have that $\varphi_d(yN)=ydy^{-1}$. By choice of $\tau$, we must have
$y\in N_G(D)$, so $ydy^{-1}\in D$. In particular,
$$\lambda(ydy^{-1})=\rho(ydy^{-1}).$$
This proves the~claim.

Therefore, $(\rho^*\circ\gamma)|_N$ and $(\lambda^*\circ\gamma)|_N$ agree
exactly on~$D$. We claim that $\rho^*\circ\gamma$ and
$\lambda^*\circ\gamma$ also agree on~$H$.

Let $h\in H$. Then $\gamma(h)=\pi(h)\varphi_h$, as in
\ref{extensions}. Since $\lambda^*$ and $\rho^*$ leave the $G/N$
component unchanged, we may concentrate on $\varphi_h$.

Let $yN\in G/N$. By definition of $\gamma$ we have
$$\eqalign{\varphi_h(yN) &= \tau\left(yN\pi(h)^{-1}\right) h \tau(yN)^{-1}\cr
&=yh'^{-1} h y^{-1}\cr}$$
where $yh'^{-1}=\tau\bigl({yN\pi(h)^{-1}}\bigr)$. In particular, we
have $h'^{-1}\equiv h^{-1} \pmod{N}$, so that $h'^{-1}h\in H\cap N$. Therefore,
$h'^{-1}h\in D$, and $y\in N_G(D)$ by choice of $\tau$, so
$\varphi_h(yN)\in D$ for all $yN\in G/N$. Therefore,
$\lambda^*\circ\varphi_h =\rho^*\circ\varphi_h$. So we have that
$\lambda^*\circ\gamma|_H =\rho^*\circ\gamma|_H$. In particular, the
dominion is contained in the equalizer of these two maps.

We therefore have that 
$$\langle H,D\rangle \subseteq {\rm dom}_G^{\cal
NQ}(H)\subseteq HN,$$
and $N\cap {\rm dom}_G^{\cal NQ}(H) = D$.  Now let $d$ be an
element of ${\rm dom}_G^{\cal NQ}(H)$. Since ${\rm dom}_G^{\cal
NQ}(H)\subseteq HN$ by \ref{mckaygeneralized}, we can write $d=hn$ for
some $h\in H$, $n\in N$. In particular, $h^{-1}d = n\in {\rm
dom}_G^{\cal NQ}(H)\cap N$. Since $D={\rm dom}_G^{{\cal NQ}}(H)\cap
N$, we must have $n\in D$, so ${\rm dom}_G^{\cal
NQ}(H)\subseteq\langle H,D\rangle$, as claimed. This proves the
equality and concludes the~proof.\endproof

\rmrk{coproductneeded} Note that we used the amalgamated coproduct
${N\amalg_D^{\cal N} N}$ in the construction. The properties that were
used of this coproduct were that there are two maps from our
group~$N$ into the amalgamated coproduct, 
which agree on~$D$ and disagree elsewhere. We could replace
the amalgamated coproduct with any group with this property. For
example, if ${\cal N}={\cal A}$, the variety of abelian groups, we
could take the group $N/D$ and the canonical projection and the zero
map from $N$ into~$N/D$. 

\rmrk{strangecond} We might wonder whether the condition that
$N_G(D)N=G$ is really necessary. It is certainly used in the
construction given in the proof. For suppose that some coset $tN\in
G/N$ does not intersect $N_G(D)$. When we try to prove that
$\lambda^*\gamma|_D = \rho^*\gamma|_D$, we run into a problem, for
there exists some $d\in D$ such that $\tau(tN)d\tau(tN)^{-1}\notin D$,
whatever the definition of $\tau(tN)$ is, and this means that
$\lambda\circ\varphi_d(tN) \not= \rho\circ\varphi_d(tN)$. In any case,
we already have examples where $HD\propcont {\rm dom}_G^{\cal
NQ}(H)$. Namely, letting ${\cal N}={\cal Q}={\cal A}$, ${\cal V}$ is
the variety of metabelian groups; since $D=H\cap {\cal Q}(G)$, the
equality given in the theorem would imply that dominions are trivial
in ${\cal A}^2$. However, there are nontrivial dominions in this
variety (see~{\bf [\cite{nildomsprelim}]}), so the equality does not
hold in~general.

\cor{bigoneprime}{Let ${\cal V}={\cal NQ}$ be a nontrivial factorization
of~${\cal V}$, and let $G\in {\cal V}$, $H$ a subgroup of~$G$, and
$N={\cal Q}(G)$. Let $D={\rm dom}_N^{\cal N}(N\cap H)$, and let $D'$
be a subgroup of $N$ such that $D\subseteq D'$, ${\rm dom}_N^{\cal
N}(D') = D'$, $HD'=D'H$ and $N_G(D')N = G$. Then
$$\langle H,D\rangle \subseteq {\rm dom}_G^{\cal NQ}(H)
\subseteq HD'.$$}

\proof Let $H'=\langle H,D'\rangle$. Since $HD'=D'H$, it follows that
$H'=HD'$. First we claim that $H'\cap N=D'$. Indeed, if $hd'\in H'\cap
N$, with $h\in H$ and $d'\in D'$, then since $D'\subseteq N$, we must
have $h\in N$. But $H\cap N=D\subseteq D'$, so $hd'\in D'$,
as~claimed.

By hypothesis, ${\rm dom}_N^{\cal N}(D')=D'$, so we apply \ref{bigone}
to~$H'$ which tells us that $H'={\rm dom}_G^{\cal NQ}(H')$. Since
$H\subseteq H'$, it follows that ${\rm dom}_G^{\cal
NQ}(H)\subseteq H'$, as~claimed\endproof

\cor{simplifone}{Let ${\cal V}={\cal NQ}$ be a nontrivial factorization
of~${\cal V}$, let $G\in {\cal V}$, $H$ a subgroup of~$G$,
$N={\cal Q}(G)$, and let $D'$ be the
normal closure of $H\cap N$ in~$G$. Then 
$$\langle H, {\rm dom}_N^{\cal N}(H\cap N)\rangle \subseteq {\rm
dom}_G^{\cal NQ}(H)\subseteq \langle H,D'\rangle.\leqno(\numbeq{sharpeq})
$$}

\proof This is just a special case of \ref{bigoneprime}, since the
normality of~$D'$ implies that $HD'=D'H$.\endproof

\cor{simpliftwo}{Let ${\cal V}={\cal NQ}$ be a nontrivial factorization
of~${\cal V}$, and let $G\in {\cal V}$, $H$ a subgroup of~$G$, and
$N={\cal Q}(G)$. If $H\cap N\triangleleft G$, then 
${\rm dom}_G^{\cal NQ}(H)=H.$ In particular, if $H\cap N=\{e\}$, then
${\rm dom}_G^{\cal NQ}(H)=H.$}

\proof If $H\cap N\triangleleft G$, then $H\cap N$ is also normal
in~$N$, hence equals its own dominion in~$N$. The rest now now follows from
\ref{bigoneprime}.\endproof

It is not hard to show that in the context of~\ref{simplifone},
either, neither, or both of the inclusions in~(\ref{sharpeq}) can
be~proper.

\Section{Applications: varieties with nontrivial
dominions}{applicationsnumber}

First, we make a few remarks for the extreme cases,
when ${\cal N}\subseteq {\cal Q}$, and when ${\cal N}\cap {\cal Q} =
{\cal E}$. In the latter case, we say that ${\cal N}$ and ${\cal Q}$
are~{\it disjoint.}

\thm{contained}{Let $\cal N$ and~$\cal Q$ be two varieties, such that
$\cal N\subseteq\cal Q$, and let $G\in\cal N$. Then for any subgroup $H$
of~$G$, ${\rm dom}_G^{\cal NQ}(H) = H$.}

\proof Follows from \ref{secondfactor}, since $G\in {\cal N}$ implies
$G\in {\cal Q}$.\endproof

Before going to the second case, let us characterize the pairs of
varieties ${\cal N}$ and~${\cal Q}$ which are disjoint. Recall that we
say that a variety $\cal N$ has {\it exponent $n$} if the variety
satisfies the law~$x^n$.

\lemma{whendisjoint}{Two varieties ${\cal N}$ and ${\cal Q}$ are
disjoint if and only if they are of relatively prime exponents.}

\proof If ${\cal N}$ and ${\cal Q}$ are of relatively prime exponents,
then any group $G\in {\cal N}\cap{\cal Q}$ satisfies the laws $x^n$
and~$x^m$, where ${\rm gcd}(n,m)=1$. Therefore, it satisfies the law
$x$, so $G=\{e\}$. Hence ${\cal N}$ and ${\cal Q}$ are disjoint.
Now suppose that ${\cal N}$ and ${\cal Q}$ are not of relatively prime
exponents. Then there must be a nontrivial factor group of the free
${\cal N}$-group of rank~1 which is isomorphic to some nontrivial
factor group of the free ${\cal Q}$-group of rank~1, so the varieties are
not disjoint.\endproof

\thm{disjointvar}{Let ${\cal V}={\cal NQ}$, where ${\cal N}$ and
{$\cal Q$} are disjoint nontrivial varieties of groups. Let
$G\in{\cal N}$, and let $H$ be a subgroup of~$G$. Then 
$${\rm dom}_G^{\cal NQ}(H) = {\rm dom}_G^{\cal N}(H).$$}

\proof Let $D={\rm dom}_{{\cal Q}(G)}^{\cal N}(H\cap {\cal
Q}(G))$. Since $G\in {\cal N}$, it follows that ${\cal Q}(G)=G$, so
$N_G(D){\cal Q}(G) = G$. The conditions of \ref{bigone} are thus
satisfied. Since ${\cal Q}(G)=G$, $H\cap {\cal Q}(G)=H$, so $D={\rm
dom}_G^{\cal N}(H)$; since $H\subseteq D$, it follows that
$$HD=D={\rm dom}_G^{\cal NQ}(H)$$
as~claimed.\endproof

Recall as well that a group~$G$ in a variety ${\cal V}$ is said to be
{\it absolutely closed} (in~${\cal V}$) if for any group $K\in{\cal
V}$, with~$G$ a subgroup of~$K$, we have
${\rm dom}_K^{\cal V}(G) = G$.

\thm{disjointvarabcl}{Let ${\cal V}={\cal NQ}$, where ${\cal N}$ and
${\cal Q}$ are disjoint nontrivial varieties of groups, and let $G\in{\cal
Q}$. Then~$G$ is absolutely closed in~${\cal V}$.}

\proof Let $K$ be a group in ${\cal V}$ such that $G$ is a subgroup
of~$K$. Consider the subgroup $G\cap {\cal Q}(K)$. Since ${\cal Q}(K)$
is an ${\cal N}$-group and $G\in {\cal Q}$, it follows that $G\cap
{\cal Q}(K)$ is trivial.
By \ref{simpliftwo}, ${\rm dom}_K^{\cal V}(G) = G$.\endproof

\thm{uncountablymany}{There are uncountably many varieties of solvable
groups with instances of nontrivial dominions. Moreover, there are
uncountably many varieties of groups with instance of
nonsurjective~epimorphisms.} 

\proof In {\bf [\cite{nildomsprelim}]} we established that the varieties of
nilpotent groups of class at most~2 and exponent $\ell^2$ (for an odd
prime~$\ell$) have instances of nontrivial dominions. Also, Ol'\v{s}anski\v{\i}
(see~{\bf [\cite{olsanskii}]}) has shown that for any two relatively
prime odd numbers $p$ and~$q$, there exist uncountably many
varieties of solvable groups with exponent $8pq$; and Vaughan-Lee
(see~{\bf [\cite{vlee}]}) has established the existence of uncountably many
varieties of solvable groups with exponent $16$. \ref{disjointvar} now gives
uncountably many varieties of solvable groups with nontrivial
dominions, by taking products of suitably chosen subvarieties
of~${\cal N}_2$ and exponent ${\ell}^2$, with solvable varieties of
exponent relatively prime to~$\ell$.

An example of B.H.~Neumann in~{\bf [\cite{pneumann}]} shows that 
${\rm Var}(A_5)$ has instances of nonsurjective epimorphisms, namely the
immersion of $A_4$ into $A_5$. Now consider the product of ${\cal
V}={\rm Var}(A_5)$ with any variety of solvable groups ${\cal S}$; by
the theorems of Ol'\v{s}anski\v{\i} and Vaughan-Lee quoted above,
there are uncountably many such varieties. Since $A_5$ is simple and
nonabelian, it follows that ${\cal S}(A_5) = A_5$ for all such ${\cal
S}$, so that \ref{bigone} applies to give
$${\rm dom}_{A_5}^{{\cal V}{\cal S}}(A_4) = {\rm dom}_{A_5}^{{\cal
 V}}(A_4) =  A_5;$$
so ${\cal V}{\cal S}$ has instances of nonsurjective epimorphisms, giving
the~result.\endproof

For the property of having nontrivial dominions, we can show even more.
Using \ref{bighalf}, we can establish
that if ${\cal N}$ is a variety with this property, and
${\cal Q}$ is any nontrivial variety, then ${\cal NQ}$ also has
instances of nontrivial dominions.

The idea is to find a group in ${\cal NQ}$ whose ${\cal Q}$-verbal
subgroup is isomorphic to a group in ${\cal N}$ which has a subgroup
with a nontrivial dominion. Then we can pick the corresponding
subgroup and apply
\ref{bighalf}. To prove this, we first need a few technical lemmas.
Along the way we will also remark on a couple of consequences of
these lemmas that are interesting in their own~right.

Some of the following lemmas are stated in generality, with ${\cal V}$
representing an arbitrary variety of algebras of a given type. In any
case, we drop for now the assumption that ${\cal V}$ can be written as
a product ${\cal V}={\cal NQ}$.

The first lemma is due to Isbell, and it is valid not only for
varieties but also for what he calls ``right closed categories'':

\lemma{fgsuffices}{{\rm (Theorem 1.2 in {\bf [\cite{isbellone}]})} Let
${\cal V}$ be a variety of algebras, and suppose
that ${\cal V}$ has instances of nontrivial dominions. Then there
exists a finitely generated algebra $G\in{\cal V}$ and a finitely 
generated subalgebra $H$ such that
$$H\propcont{\rm dom}_{G}^{\cal V}(H).$$}

\proof We sketch the argument, and refer the reader to Isbell's paper
for details. Let $G'$ be an algebra, and $H'$ a subalgebra such that
the dominion of~$H'$ in~$G'$ is nontrivial; let $d\in G'\setminus H'$ be
an element of the dominion not in~$H'$. One considers the coproduct
$G'\amalg^{\cal V} G'$, the amalgamated coproduct $G'\amalg^{\cal
V}_{H'} G'$, and the natural quotient map $\pi$ between them. The fact
that $d$ is in the dominion means that the elements $\lambda(d)$ and $\rho(d)$ 
of $G'\amalg^{\cal V} G'$ have the same image under~$\pi$.
This can be expressed as a finite sequence of elements
$(w_0,\ldots,w_n)$ of $G'\amalg^{\cal V}G'$, starting with
$w_0=\lambda(d)$, ending with $w_n=\rho(d)$, and such that in 
$(G'\amalg^{\cal V} G')\times (G'\amalg^{\cal V}G')$, each element
$(w_j,w_{j+1})$ 
lies in the subalgebra generated by all elements of the form $(x,x)$,
$(\lambda(h),\rho(h))$, and $(\rho(h),\lambda(h))$, with $h\in H'$. 
We let $H$ be the algebra
generated by the $h$'s used in this finite expression, and let $G$ be
obtained by throwing in all other elements that we need.\endproof 

\rmrk{finitebases} Note also that at most finitely many identities
of~${\cal V}$ are needed to express the fact that $\lambda(d)$ and
$\rho(d)$ have the same image under~$\pi$.  So the argument above
shows that if the variety ${\cal V}$ is not finitely based (that is,
it is not defined by a finite number of identities), then there exists
a finitely based variety ${\cal V}'$, with ${\cal V}\subseteq {\cal
V}'$, and a finitely generated algebra $G\in {\cal V}$ with a finitely
generated subalgebra $H$, such that the dominion of~$H$ in~$G$ (in the
variety ${\cal V}'$) is nontrivial. Namely, we let ${\cal V}'$ be the
variety defined only by the identities used in the process described.

\rmrk{finitepresentation} The
argument may also be adapted to show
that one can take~$H$ finitely generated, and $G$ not merely finitely
generated, but also finitely~presented.

\rmrk{notallcats} The statement analogous to \ref{fgsuffices}  is {\it
not} true if we replace the variety ${\cal V}$ with an 
arbitrary category of algebras, and not even in the case of
pseudovarieties. Explicitly, Example~8.87 in~{\bf
[\cite{nildomsprelim}]} shows that there are nontrivial dominions in
the category of all nilpotent groups (which is a
pseudovariety). However, Theorem~3.11 in~{\bf [\cite{fgnilprelim}]}
shows that dominions of subgroups of finitely generated nilpotent
groups are trivial in the category of all nilpotent~groups.  Thus
\ref{fgsuffices} does not hold in general.

Recall that the variety ${\cal V}^{(m)}$ is the variety defined by the
$m$-variable identities of~${\cal V}$, or equivalently, the variety of
all algebras whose $m$-generator algebras belong to~${\cal V}$.

\lemma{suffhighn}{Let ${\cal V}$ be a variety of algebras with instances
of nontrivial dominions. Then there exists an $n>1$ such that for all
$m\geq 2n$, the variety~${\cal V}^{(m)}$ also has instances of
nontrivial~dominions.}

\proof By \ref{fgsuffices}, there is a finitely generated algebra~$G$
with a subalgebra~$H$ which has nontrivial dominion. Say that $G$ is
an $n$-generated algebra. Let $m\geq 2n$, and consider the algebra
$$K = G \amalg^{{\cal V}^{(m)}} G.$$

Since $K$ is generated by two copies of~$G$, it is a $2n$-generated
algebra, hence lies in~${\cal V}$. In particular, since the two
embeddings of $G$ into~$K$ agree on~$H$, they agree on ${\rm
dom}_G^{\cal V}(H)$, which properly contains~$H$. Hence $H$ is
properly contained in ${\rm dom}_G^{{\cal V}^{(m)}}(H)$,
as~claimed.\endproof

\rmrk{maybenotnew} We remark, however, that we may have
${\cal V}^{(m)}={\cal V}$ for sufficiently large~$m$ (for example, if
${\cal V}$ is finitely based), so \ref{suffhighn} may not
yield any new~information.

We now return our attention to varieties of groups.

\lemma{freesuffices}{Let ${\cal V}$ be a variety of groups with
instances of nontrivial dominions. Then there exists an $n>1$, and a
subgroup $H$ of $F_n=F_n({\cal V})$, the relatively free
${\cal V}$-group of rank~$n$, such~that
$H\propcont {\rm dom}_{F_n}^{\cal V}(H)$.}

\proof By \ref{fgsuffices}, there exists a finitely generated group
$G\in {\cal V}$, and a subgroup $K$ of~$G$, such that
$K\propcont {\rm dom}_G^{\cal V}(H)$.

Let $g_1,\ldots,g_n$ be generators for $G$, and let $\pi\colon
F_n({\cal V})\to G$ be the map from the relatively free group of
rank~$n$ to~$G$ sending the free generators of $F_n({\cal V})$,
$x_1,\ldots,x_n$ to $g_1,\ldots,g_n$, respectively.
Let $N={\rm ker}(\pi)$. 

Let $H$ be the subgroup of~$F_n({\cal V})$ containing~$N$ which
corresponds, under $\pi$, to the subgroup $K$ of~$G$. Since dominions
respect quotients in varieties,
$${\rm dom}_G^{\cal V}(K) = {\rm dom}_{F_n/N}^{\cal V}(H/N)
= \bigl({\rm dom}_{F_n}^{\cal V}(H)\bigr) / N.$$

Since $K\propcont {\rm dom}_G^{\cal V}(K)$, it follows that
$H\propcont {\rm dom}_{F_n}^{\cal V}(H)$, as desired.\endproof

\rmrk{ingeneralidont} The proof of \ref{freesuffices} does not easily
generalize to an arbitrary variety of algebras; the main problem lies
in the argument regarding the quotient. It is easy to see that, given
an algebra $A$, subalgebra $B$, and congruence relation~$R$, we have
an inclusion
$${\rm dom}_{G}^{\cal V}(H)\bigm/R \subseteq {\rm dom}_{G/R}^{\cal
V}(H/R);$$ in the case of varieties of groups, the reverse inclusion
also holds, but this may not be true in general (it may not even be
true of arbitrary classes of groups, although it is true of unions of
pseudovarieties {\bf [\cite{nildomsprelim}]}). In the argument
above, we know that the subalgebra on the right hand side is bigger
than $H/R$, but if we simply pull back to~$G$ we could, conceivably,
end up with a subalgebra which is equal to $H$. I~do not know if the
general result itself is true, using a different proof.

\lemma{finftysuffices}{If ${\cal V}$ is a variety of groups with
instances of nontrivial dominions, then there exists a subgroup $H$ of
the countably generated relatively free ${\cal V}$-group
$F=F_\infty({\cal V})$ such that
$H\propcont {\rm dom}_{F}^{\cal V}(H)$.}

\proof By \ref{freesuffices}, there exists an $n>1$, and a subgroup
$H$ of $F_n=F_n({\cal V})$ such that $H\propcont {\rm dom}_{F_n}^{\cal
V}(H)$. Let $i\colon F_n({\cal V})\to F$ be the immersion of $F_n$
into $F$ sending the free generators of $F_n$ to the first $n$ free
generators of~$F$. Then
$$i(H)\propcont {\rm dom}_{i(F_n)}^{\cal V}(i(H)) \subseteq
{\rm dom}_{F}^{\cal V}(i(H)).$$
This proves the lemma.\endproof

\lemma{verbalofinfty}{Let $F=F_{\infty}({\cal G})$ be the absolutely
free group on infinitely many generators $x_1,x_2,\ldots,x_n,\ldots$, and let
${\bf V}$ be a nontrivial set of words. Then ${\bf V}(F)$ is free of
countably infinite~rank.}

\proof ${\bf V}(F)$ is a subgroup of an absolutely free group, so by a
classical theorem of Schreier it is itself free. It remains to show that
it is not finitely generated.

Suppose that ${\bf V}(F)$ is finitely generated. Then there is a bound
to the subscripts of the free generators of~$F$ that occur in the generators
of~${\bf V}(F)$, so there exists an $n_0>0$ such that all generators
of~$V(F)$ involve only $x_1$, $x_2,\ldots,x_{n_0}$. 

Now let ${\bf v}(z_1,\ldots,z_r)\in {\bf V}$ be a nonidentity element; note
that the element $${\bf v}(x_{n_0+1},\ldots,x_{n_0+r})$$
is a
nonidentity element of~${\bf V}(F)$, expressible in terms of a set of
generators disjoint from $x_1,\ldots,x_{n_0}$, which contradicts the
assumption that all elements of~${\bf V}(F)$ are expressible in terms
of $x_1,\ldots,x_{n_0}$ only. We conclude that ${\bf V}(F)$ is of
countably infinite rank, as~claimed.\endproof

We can now prove the result we mentioned above:

\thm{bigtwo}{Let ${\cal N}$ be a variety of groups with instances of
nontrivial dominions, and let ${\cal Q}$ be any variety of groups
different from ${\cal G}$. Then ${\cal NQ}$ also has instances of
nontrivial~dominions.}

\proof Let $F=F_{\infty}({\cal NQ})$ be the countably generated relatively
free group in~${\cal NQ}$. By \ref{verbproduct} $F$ is given by
$F_{\infty}({\cal G})/{\cal N}({\cal Q}(F_{\infty}({\cal G})))$, and
the ${\cal Q}$-verbal subgroup of $F$ is given by 
$${\cal Q}(F) = {\cal Q}(F_{\infty}({\cal G}))/{\cal N}({\cal
Q}(F_{\infty}({\cal G}))).$$

Since ${\cal Q}\not={\cal G}$, we have ${\cal
Q}(F_{\infty})\not=\{e\}$. We claim that ${\cal Q}(F)\cong
F_{\infty}({\cal N})$. Indeed, ${\cal
Q}(F_{\infty}({\cal G}))$ is isomorphic to $F_{\infty}({\cal G})$ by
\ref{verbalofinfty}. Therefore,
$$\eqalign{{\cal Q}(F) &= {\cal Q}(F_{\infty}({\cal G}))/ {\cal N}
({\cal Q}(F_{\infty}({\cal G})))\cr
&\cong F_{\infty}({\cal G})/{\cal N}(F_{\infty}({\cal G}))\cr
&\cong F_{\infty}({\cal N}).\cr}$$

By \ref{finftysuffices}, since ${\cal N}$ has instances of nontrivial
dominions, there exists a subgroup $H$ of $F_{\infty}({\cal N})$ such
that
$$H\propcont {\rm dom}_{F_{\infty}({\cal N})}^{\cal N}(H).$$
Now consider $H$ as a subgroup of ${\cal Q}(F)$, and hence as a
subgroup of~$F$. By \ref{bighalf}, if we let
$G=F$, we have that
$$\Bigl\langle H, {\rm dom}_{{\cal Q}(F)}^{\cal N}\bigl(H\cap {\cal
Q}(F)\bigr)\Bigr\rangle \subseteq {\rm dom}_{F}^{\cal NQ}(H).$$

But $H\cap {\cal Q}(F) = H$, and by choice of~$H$, $H\propcont {\rm
dom}_{{\cal Q}(F)}^{\cal N}(H)$. In particular,
$$H\propcont \Bigl\langle H, {\rm dom}_{{\cal Q}(F)}^{\cal
N}\bigl(H)\Bigr\rangle = {\rm dom}_{{\cal Q}(F)}^{\cal N}(H) \subseteq
{\rm dom}_{F}^{\cal NQ}(H)$$ 
so ${\cal NQ}$ has instances of nontrivial dominions, as~claimed.\endproof

\rmrk{secondproof} \ref{bigtwo} now provides an alternative proof to
the first sentence of \ref{uncountablymany}, simply by taking any
solvable variety ${\cal N}$ with instances of nontrivial dominions
(e{.}g{.} the variety of nilpotent groups of class 2, ${\cal A}^2$,
etc{.}), and looking at all possible products ${\cal NQ}$ with ${\cal
Q}$ solvable and nontrivial.

\rmrk{converse} We might ask whether there is a converse to
\ref{bigtwo}. That is, given a nontrivial factorization ${\cal V}=
{\cal NQ}$, can we deduce that ${\cal V}$ has trivial dominions
if~${\cal N}$ (or maybe, if both ${\cal N}$ and~${\cal Q}$) has trivial
dominions? The answer to this question is no, as we have already
noted that there are nontrivial dominions in the variety of
metabelian groups, which has nontrivial factorization ${\cal V}={\cal
A}{\cal A}$, and all dominions are trivial in the variety of abelian~groups.

\cor{solvable}{Let $\ell\geq 2$, and let ${\cal S}_{\ell}={\cal
A}^{\ell}$ be the variety of all solvable groups of solvability length
at most $\ell$. Then ${\cal S}_{\ell}$ has instances of nontrivial
dominions.}

\proof We write ${\cal S}_{\ell} = {\cal A}^2{\cal A}^{\ell -2}$. As
we have already mentioned, there are instances of nontrivial dominions
in~${\cal A}^2$. The result now follows from \ref{bigtwo}.\endproof

\rmrk{allsolvable} On the other hand, if ${\cal C}_{\cal S}$ is the
category of all solvable groups, then dominions are trivial in ${\cal
C}_{\cal S}$. For given $G\in {\cal C}_{\cal S}$, and $H$ a subgroup
of~$G$, we know that $G\in {\cal S}_{\ell}$ for some $\ell>0$.
Hence $G\in{\cal S}_{\ell+1} = {\cal S}{\cal S}_{\ell}$. By
\ref{secondfactor}, ${\rm dom}_G^{{\cal S}_{\ell+1}}(H)=H$. In
particular, ${\rm dom}_G^{{\cal C}_{\cal S}}(H)=H$. 

\ref{solvable} shows the converse of
\ref{bigtwo} is false. This suggests 
the following questions: Given a nontrivial variety ${\cal N}$, does
${\cal N}^2$ always have instances of nontrivial dominions? More
generally, given a variety ${\cal V}={\cal NQ}$, with ${\cal N}, {\cal
Q}$ nontrivial, will ${\cal V}$ necessarily have instances of
nontrivial dominions? I do not know the answer to these two questions.

%
\ifnum0<\citations{\par\bigbreak
\filbreak{\bf References}\par\frenchspacing}\fi
%
%
\ifundefined{xbirkhoff}\else
\item{\bf [\refer{birkhoff}]}{Garrett Birkhoff, {On the structure
of abstract algebras.} {\it Proc.\ Cambridge\ Philos.\ Soc.} {\bf
31} (1935), \hbox{433--454}.}\par\filbreak\fi
\ifundefined{xisbellone}\else
\item{\bf [\refer{isbellone}]}{J. R. Isbell, {Epimorphisms and
dominions.} In {\it Proc.~of the Conference on Categorical Algebra, La
Jolla 1965,} 
(Lange and Springer, New
York~1966). MR:35\#105a (The statement of the
Zigzag Lemma for {\it rings} in this paper is incorrect. The correct
version is stated in~{\bf [\cite{isbellfour}]}.)}\par\filbreak\fi
\ifundefined{xisbellfour}\else
\item{\bf [\refer{isbellfour}]}{J. R. Isbell, {Epimorphisms and
dominions IV,} {\it J.\ London Math.\ Soc.~(2),}
{\bf 1} (1969) \hbox{265--273.} {MR:41\#1774}}\par\filbreak\fi
\ifundefined{xwreathext}\else
\item{\bf [\refer{wreathext}]}{L.~Kaloujnine and Marc Krasner,
{Produit complet des groupes de permutations et le probl\`eme
d'extension des groupes III,} {\it Acta Sci.\ Math.\ Szeged} {\bf 14}
(1951) \hbox{69--82}. {MR:14,242d}}\par\filbreak\fi
\ifundefined{xmckay}\else
\item{\bf [\refer{mckay}]}{Susan McKay, {Surjective epimorphisms
in classes
of groups,} {\it Quart.\ J.\ Math.\ Oxford (2),\/} {\bf 20} (1969),
\hbox{87--90.} {MR:39\#1558}}\par\filbreak\fi
\ifundefined{xnildomsprelim}\else
\item{\bf [\refer{nildomsprelim}]}{Arturo Magidin, {Dominions in
varieties of nilpotent groups,} {\it Comm.\ Alg.} to appear.}\par\filbreak\fi
\ifundefined{xfgnilprelim}\else
\item{\bf [\refer{fgnilprelim}]}{Arturo Magidin, {Dominions in
finitely generated nilpotent groups,} {\it Comm.\
Alg.} to appear.}\par\filbreak\fi
\ifundefined{xhneumann}\else
\item{\bf [\refer{hneumann}]}{Hanna Neumann, {\it Varieties of
Groups,} {(Ergebnisse der Mathematik und ihrer Grenz\-ge\-biete,\/
New series, Vol.~37, Springer Verlag~1967)}. {MR:35\#6734}}\par\filbreak\fi
\ifundefined{xpneumann}\else
\item{\bf [\refer{pneumann}]}{Peter M.~Neumann, {Splitting groups
and projectives
in varieties of groups,} {\it Quart.\ J.\ Math.\ Oxford} (2), {\bf
18} (1967),
\hbox{325--332.} {MR:36\#3859}}\par\filbreak\fi
\ifundefined{xolsanskii}\else
\item{\bf [\refer{olsanskii}]}{A. Ju.~Ol'\v{s}anski\v{\i}, {On the
problem of a finite basis of identities in groups,} {\it
Izv.\ Akad.\ Nauk.\ SSSR} {\bf 4} (1970) no. 2
\hbox{381--389.}}\par\filbreak\fi
\ifundefined{xvlee}\else
\item{\bf [\refer{vlee}]}{M. R. Vaughan-Lee, {Uncountably many
varieties of groups,} {\it Bull.\ London Math.\ Soc.} {\bf 2} (1970)
\hbox{280--286.} {MR:43\#2054}}\par\filbreak\fi
\ifnum0<\citations\nonfrenchspacing\fi

\bigskip
 
\vfill\eject
\immediate\closeout\aux
\end